\documentclass[a4paper,11pt,reqno,noindent]{amsart}
\usepackage[centertags]{amsmath}
\usepackage{amsfonts,amssymb,amsthm,amscd,dsfont,cases,esint,enumerate,xcolor}
\usepackage[T1]{fontenc}
\usepackage[english]{babel}
\usepackage[applemac]{inputenc}
\usepackage[body={15cm,21.5cm},centering]{geometry} 
\usepackage{fancyhdr}
\numberwithin{equation}{section}

\theoremstyle{plain}
\newtheorem{theor10}{Theorem}

\newtheorem{prop10}[theor10]{Proposition}

\newtheorem{cor10}[theor10]{Corollary}

\newtheorem{lem10}[theor10]{Lemma}

\newtheorem{theor0}{Theorem}[section]
\newenvironment{theor}
  {\pushQED{\qed}\begin{theor0}}
  {\popQED\end{theor0}}
\newtheorem{lem0}[theor0]{Lemma}

\newtheorem{prop0}[theor0]{Proposition}

\newtheorem{cor0}[theor0]{Corollary}

\theoremstyle{definition}
\newtheorem{rems0}[theor0]{Remarks}

\newtheorem{rem0}[theor0]{Remark}

\theoremstyle{plain}
\newtheorem{as0}[theor0]{Assumption}

\newtheorem*{asn0*}{\assumptionnumber}
  \providecommand{\assumptionnumber}{}
\makeatletter
\newenvironment{asn0}[2]
   {\renewcommand{\assumptionnumber}{Assumption \!#1 {\normalfont--- #2}}
    \begin{asn0*}
    \protected@edef\@currentlabel{{\normalfont#1}}}
   {\end{asn0*}}
\makeatother

\makeatletter
\newenvironment{asn01}[1]
   {\renewcommand{\assumptionnumber}{Assumption \!#1}
    \begin{asn0*}
    \protected@edef\@currentlabel{{\normalfont#1}}}
   {\end{asn0*}}
\makeatother

\newcommand{\Pc}{\mathcal{P}}

\newcommand{\Bc}{\mathcal{B}}
\newcommand{\Sc}{\mathcal S}

\newcommand{\R}{\mathbb R}

\newcommand{\loc}{{\operatorname{loc}}}

\newcommand{\Id}{\operatorname{Id}}

\newcommand{\Ld}{\operatorname{L}}

\newcommand{\step}[1]{\noindent \textit{Step} #1.}

\newcommand{\Pm}{\mathbb{P}}
\newcommand{\pr}[1]{\mathbb{P}\left[ #1 \right]}

\newcommand{\expec}[1]{\mathbb{E}\left[ #1 \right]}

\newcommand{\expecm}[1]{\mathbb{E}\big[ #1 \big]}
\newcommand{\expecM}[1]{\mathbb{E}\bigg[ #1 \bigg]}

\usepackage[colorlinks,citecolor=black,urlcolor=black]{hyperref}

\title{The Clausius--Mossotti formula}

\author[M. Duerinckx]{Mitia Duerinckx}
\address[Mitia Duerinckx]{Universit\'e Libre de Bruxelles, D\'epartement de Math\'ematique, 1050~Brussels, Belgium}
\email{mitia.duerinckx@ulb.be}
\author[A. Gloria]{Antoine Gloria}
\address[Antoine Gloria]{Sorbonne Universit\'e, CNRS, Universit\'e de Paris, Laboratoire Jacques-Louis Lions, 75005~Paris, France \& Institut Universitaire de France \& Universit\'e Libre de Bruxelles, D\'epartement de Math\'ematique, 1050~Brussels, Belgium}
\email{antoine.gloria@sorbonne-universite.fr}

\begin{document}
\selectlanguage{english}

\maketitle

\begin{abstract}
In this note, we provide a short and robust proof of the Clausius--Mossotti formula for the effective conductivity in the dilute regime, together with an optimal error estimate. The proof makes no assumption on the underlying point process besides stationarity and ergodicity, and can be applied to dilute systems in many other contexts.
\end{abstract}

\maketitle

\section{Effective conductivity problem}
We start by recalling the notion of effective conductivity in the sense of stochastic homogenization theory for an heterogeneous material made of inclusions in a given matrix of homogeneous conductivity. Let $d\ge1$ denote the space dimension.

\subsection{Stochastic setting}
We shall use a statistical description for the set of inclusions in the material. Restricting to spherical inclusions for notational simplicity, we let
\begin{equation}\label{eq:inclusion}
\Bc(\Pc)\,:=\,\cup_{x\in\Pc}B(x),
\end{equation}
where $B(x)=B+x$ stands for the unit ball centered at $x$ in $\R^d$, and where $\Pc$ is the set of centers of the inclusions. It remains to define statistical ensembles for the latter.
We call {\it point set} any countable subset $\Pc\subset\R^d$ that is {locally finite} in the sense that for any bounded $E\subset\R^d$ the number of points of~$\Pc$ in~$E$ is finite, $\Pc(E):=\sharp\{\Pc \cap E\}<\infty$.
A point set~$\Pc$ can be represented by the associated locally finite measure $\sum_{x \in \Pc} \delta_x$, which acts on the space of compactly supported continuous functions via $f \mapsto\Pc(f):=\sum_{x\in \Pc} f(x)$. We endow the space $\Omega$ of point sets with the smallest $\sigma$-algebra that makes all evaluation maps $\Pc\mapsto f(\Pc)$ measurable.
A {\it random point process} is then defined as a probability measure $\Pm$ on $\Omega$, and we denote by $\expec{\cdot}$ the associated expectation.
We further define stationarity and ergodicity with respect to translations $\Pc+z:=\{x+z:x\in\Pc\}$ of point sets: The point process is said to be {\it stationary} (or {\it statistically translation-invariant}) if for any measurable set $A\subset\Omega$ we have $\pr{A+z}=\pr{A}$ for all $z\in\R^d$, where we use the notation $A+z:=\{\Pc+z:\Pc\in A\}$. The process is said to be {\it ergodic} if any measurable set $A\subset\Omega$ that is translation-invariant, in the sense that $\pr{A\setminus(A+z)}=\pr{(A+z)\setminus A}=0$ for all $z\in\R^d$, satisfies~$\pr{A}=0$ or~$1$.
In the sequel, the set of inclusions in the material is modeled by~\eqref{eq:inclusion} with $\Pc$ sampled according to some stationary and ergodic random point process. Note that the inclusions are allowed to overlap in general.

\subsection{Effective conductivity}
Given a stationary ergodic random point process as defined above, we consider the associated coefficient field
\begin{equation}\label{e.coeff}
A(y):=A_1 +\mathds{1}_{\Bc(\Pc)}(y)(A_2-A_1) ,\qquad y\in\R^d,
\end{equation}
where $A_1,A_2\in\R^{d\times d}$ are two strongly elliptic matrices and where $\Pc$ is sampled according to the point process.
This models a homogeneous material of conductivity $A_1$ that is perturbed by disordered spherical inclusions of another material of conductivity $A_2$.
On large scales, in the sense of the homogenization theory, this two-phase heterogeneous material behaves like a homogeneous material with some effective conductivity $\bar A$ defined  by
\begin{equation}\label{0-1}
\bar A e = \expec{A(\nabla \phi_e+e)},
\end{equation}
where $\phi_e$ is the so-called corrector defined as the unique weak solution of the whole-space equation 
\begin{equation}\label{0-2}
-\nabla \cdot A (\nabla \phi_e+e)=0,\qquad \text{in $\R^d$},
\end{equation}
in the following class: $\phi_e$ is almost surely in $H^1_\loc(\R^d)$, satisfies the anchoring condition~$\fint_{B(0)} \phi_e=0$, and its gradient $\nabla \phi_e$ is a stationary random field with vanishing expectation $\expec{\nabla \phi_e}=0$, and finite second moments~$\expec{|\nabla \phi_e|^2}<\infty$; see \cite{MR542557,MR712714,MR1329546}.
Note that $\bar A$ is not explicit in general.

\section{Dilute homogenization}

\subsection{The Clausius--Mossotti formula}
Homogenization was neither born in the mathematical community in the 1970s, nor in the engineering community the decade before: it emerged much earlier, in the second half of the 19th century, in the physics community, in the context of two-phase dispersed media.
Motivated by the works of Poisson~\cite{Poisson} and Faraday~\cite{Faraday}, Mossotti and Clausius were the first to investigate the question of the effective dielectric constant of a homogeneous background material perturbed by sparse spherical inclusions~\cite{Mossotti-36,Mossotti-50,Clausius-79}.
The problem was largely
revisited by Maxwell~\cite{Maxwell-81} for the effective conductivity of two-phase media; we refer to~\cite{Markov-00} for a detailed account of the historical context.
In modern language, these authors argued that the effective conductivity~\eqref{0-1} associated with the two-phase model~\eqref{e.coeff} takes on the following guise, in dimension $d\ge1$, in case of isotropic conductivities $A_1=\alpha\Id$ and $A_2=\beta\Id$,
\begin{equation}\label{0-10}
\bar A=\alpha \Id +\varphi \frac{\alpha d(\beta-\alpha)}{\beta+\alpha(d-1)} \Id+o(\varphi),\qquad\text{as $\varphi\downarrow0$},
\end{equation} 
where $\varphi$ stands for the volume fraction of the inclusions, that is, by the ergodic theorem,
\begin{equation}\label{eq:def-varphi}
\varphi\,:=\,\expecm{\mathds1_{\Bc(\Pc)}}\,=\,\lim_{R\uparrow\infty}\tfrac{|\Bc(\Pc)\cap RB|}{|RB|},\qquad\text{for $\Pm$-almost all $\Pc$}.
\end{equation}
In case of disjoint inclusions, note that $\varphi=\lambda|B|$ where $\lambda$ is the intensity of the point process (defined by $\lambda|E|=\expec{\Pc(E)}$ for any Borel set $E\subset\R^d$).
This dilute approximation~\eqref{0-10} for the effective medium is known as the Clausius--Mossotti formula and was soon adapted to various other physical settings, in particular by Lorenz and Lorentz for the effective refractive index of two-phase media in optics~\cite{Lorenz-80,Lorentz-09}, and by Einstein for the effective viscosity of a Stokes fluid with a dilute suspension of rigid particles~\cite{Einstein-05,Einstein-11}. Einstein's result was actually part of his PhD thesis,
where he used it to design a celebrated experiment to measure the Avogadro number; see e.g.~the inspiring historical account in~\cite{Straumann-05}.

We start by describing a heuristic argument for~\eqref{0-10}.
Recalling that~$\expec{\nabla \phi_e}=0$,
the effective conductivity~\eqref{0-1} can be decomposed as
\begin{eqnarray*}
e \cdot \bar A e &=& \expec{e\cdot A(\nabla \phi_e+e)}\\
&=&\expec{e\cdot \big(\alpha+(\beta-\alpha)\mathds 1_{\Bc(\Pc)}\big) (\nabla \phi_e+e)}\\
&=& \alpha|e|^2+ (\beta-\alpha)\,\expec{\mathds 1_{\Bc(\Pc)}\, e\cdot (\nabla \phi_e+e)},
\end{eqnarray*}
and thus, by the ergodic theorem, assuming for simplicity that inclusions are almost surely disjoint,
\begin{equation}\label{0-12}
e \cdot \bar A e
\,=\,\alpha|e|^2+ (\beta-\alpha)\,\lim_{R\uparrow \infty} \tfrac1{|RB|} \sum_{x \in \Pc \cap RB}e \cdot\int_{B(x)} (  \nabla \phi_e+e).
\end{equation}
In the dilute regime~\mbox{$\varphi\ll 1$}, inclusions are typically far from one another and therefore do not `interact' much when solving the corrector equation~\eqref{0-2}.
More precisely, for all~$x\in\Pc$, we may heuristically approximate the corrector in the inclusion $B(x)$ by the solution of a corresponding single-inclusion problem,
\begin{equation}\label{0-13}
\nabla \phi_e|_{B(x)} \simeq \nabla \psi_e(\cdot-x),
\end{equation}
where $\psi_e$ is the unique weak solution in $\dot H^1(\R^d)$ of the whole-space single-inclusion equation
\begin{equation*}
-\nabla \cdot \big(\alpha+ (\beta-\alpha)  \mathds1_{B(0)}\big) (\nabla \psi_e+e)=0,\qquad\text{in $\R^d$}.
\end{equation*}
In the present case of spherical inclusions and isotropic conductivity, this equation is explicitly solvable in form of
\begin{equation}\label{0-11}
\nabla \psi_e(x)=
\begin{cases}
Ke,&\text{for $|x|<1$};\\
\frac{K}{|x|^d}\left(e-d\frac{x\cdot e}{|x|}\frac{x}{|x|}\right)
,&\text{for $|x|>1$};
\end{cases}
\end{equation}
with 
\begin{equation*} 
K\,=\,\frac{\alpha-\beta}{\beta+\alpha(d-1)}.
\end{equation*}
Inserting this form into~\eqref{0-13} and~\eqref{0-12}, the Clausius--Mossotti formula~\eqref{0-10} heuristically follows.

\subsection{Main result}
The aim of the present note is to prove the Clausius--Mossotti formula~\eqref{0-10} in the most general setting possible and to establish a sharp error bound.
The above heuristic argument indicates that the error mostly comes from `interactions' between inclusions, as it amounts to locally replacing the corrector by solutions of single-inclusion problems, cf.~\eqref{0-13}.
To quantify this error, we need to recall the notion of {\it second-order intensity} $\lambda_2$ of the point process, which we introduced in~\cite{2008.03837}. For that purpose, we first define the minimal lengthscale $\ell\ge0$ of the point process,
\begin{equation}\label{0-14}
\ell\,:=\, \inf_{x,y\in\Pc\atop x\ne y} |x-y|_\infty,\qquad\text{for $\Pm$-almost all $\Pc$,}
\end{equation}
which is deterministic by ergodicity.
In case $\ell>0$ (that is, if the point process is hardcore), the {\it second-order intensity} is defined as
\begin{equation}\label{0-15}
\qquad\lambda_2\,:=\,\sup_{z_1,z_2\in\R^d}\,\expecM{\sum_{x_1,x_2\in\Pc\atop x_1\ne x_2}\ell^{-d}\,\mathds1_{Q_{\ell}(z_1)}(x_1)\,\ell^{-d}\,\mathds1_{Q_{\ell}(z_2)}(x_2)},
\end{equation}
where $Q_r(z):=r Q+z$ is the cube of sidelength $r$ centered at $z$. Note that, by definition~\eqref{0-14}, each cube $Q_\ell(z)$ contains at most one point of $\Pc$ almost surely.
In other words, $\lambda_2$ is the maximum expected number of couples of points that lie in the $\ell$-neighborhood of a given element of $(\R^d)^2$, properly normalized by $\ell$.
Alternatively, recalling that the {\it$2$-point density} is the non-negative function~$f_{2}$ defined 
by the following relation,
\begin{equation*}
\qquad\expecM{\sum_{x_1, x_2\in\Pc\atop x_1\ne x_2}\zeta(x_{1},x_{2})} \,=\, \int_{(\R^d)^2}\zeta f_2,\qquad\text{for all $\zeta \in C_c((\R^d)^2)$},
\end{equation*}
the definition~\eqref{0-15} of $2$-point intensity can be reformulated as
\begin{equation}\label{eq:high-intens/Re}
\lambda_2\,=\,\sup_{z_1,z_2\in\R^d}\,\fint_{Q_{\ell}(z_1)\times Q_\ell(z_2)}f_2.
\end{equation}
In case $\ell=0$ (that is, if the point process is not hardcore), this definition is naturally extended to $\lambda_2:=\|f_2\|_{\Ld^\infty((\R^d)^2)}$. For a Poisson point process, due to the tensor structure, the $2$-point intensity is simply the square of the intensity, $\lambda_2=\lambda^2$, but for a general mixing point process it can be anything in the interval $[\lambda^2,\lambda]$ and its smallness describes some form of local independence.

We can now state the main result of this note, which is an adaptation of our recent work on Einstein's formula~\cite[Theorem~1]{2008.03837} to the effective conductivity problem.
It proves the validity of the Clausius--Mossotti formula~\eqref{0-10} for the first time in the setting of general point processes.
The error bound $O(\lambda_2|\!\log\lambda_2|)$ below is new and is sharp in general, cf.~\cite[Theorem~7]{2008.03837}, showing that the approximation~\eqref{0-10} is only valid in the dilute regime if~$\lambda_2|\!\log\lambda_2|=o(\varphi)$.
In the specific case of a Poisson point process, this result follows from~\cite{CRMATH_2022__360_G8_909_0} with the improved error bound $O(\lambda_2)=O(\lambda^2)$ (without logarithmic correction). The particular case of dilute point processes obtained by Bernoulli deletion or by dilation of a given process was already treated in~\cite{MR3458165,Pertinand}.
Note that inclusions here are allowed to overlap and that no upper bound is assumed on the number of points per unit volume.

\begin{theor}\label{th:CM-DG}
Given a stationary and ergodic point process, and given $\alpha,\beta>0$,
the effective conductivity~\eqref{0-1} associated with the two-phase model~\eqref{e.coeff} with isotropic conductivities $A_1=\alpha\Id$ and $A_2=\beta\Id$ satisfies the following quantitative version of~\eqref{0-10},
\begin{equation*}
\Big|\bar A-\Big(\alpha \Id + \varphi\frac{\alpha d(\beta-\alpha)}{\beta+\alpha(d-1)} \Id \Big)\Big| \,\lesssim\,\lambda_2|\!\log\lambda_2|.
\qedhere
\end{equation*}
\end{theor}

We shall prove this result in the following slightly more general form, where conductivities $A_1$ and $A_2$ are no longer assumed to be isotropic and where $A_2$ may itself be heterogeneous.
We consider spherical inclusions for notational convenience, but we emphasize that, as in~\cite{MR3458165,2008.03837}, this is not essential (only the above explicit form of the Clausius--Mossotti formula then needs to be changed).

\begin{theor} \label{th:CM-DG+}
Let $A_1\in\R^{d\times d}$ be a strictly elliptic (non-necessarily symmetric) matrix, and let $A_2$ be a
stationary and ergodic random field of uniformly elliptic (non-necessarily symmetric) matrices.
Given a stationary and ergodic point process that is independent of~$A_2$, consider the associated coefficient field
\begin{equation}\label{eq:model-A}
A(y):=A_1+\mathds{1}_{\Bc(\Pc)}(y)(A_2(y)-A_1),\qquad y\in\R^d,
\end{equation}
where $\Pc$ is sampled according to the point process.
Then, the effective coefficient $\bar A$ associated via~\eqref{0-1} satisfies the following expansion,
\begin{align} \label{0-26+}
&\big|\bar A-(A_1 +\varphi\widehat A_2 )\big| \,\lesssim\,\lambda_2 |\!\log\lambda_2|,
\end{align}
where the first-order effective correction $\widehat A_2$ is given by
\begin{equation}\label{eq:def-corrA2}
\widehat A_2 e \,:=\,\expecM{\fint_{B} (A_2-A_1)(\nabla \psi_e+e)},
\end{equation}
where $\psi_e$ is the unique weak solution in $\dot H^1(\R^d)$ of the whole-space single-inclusion problem
\begin{equation}\label{eq:single-incl-1}
-\nabla \cdot \big(A_1+\mathds 1_{B}(A_2-A_1)\big)(\nabla \psi_e+e)=0,\qquad\text{in $\R^d$.}\qedhere
\end{equation}
\end{theor}

\subsection{Previous contributions}

The asymptotic analysis of the effective conductivity in case of a periodic array of inclusions with a small volume fraction was first addressed by Berdichevski{\u\i}~\cite{Berdichevski-83}; see also~\cite[Section~1.7]{MR1329546}.
The first justification of the  Clausius--Mossotti formula in a random setting is due to Almog in dimension $d=3$, whose results in~\cite{Almog-13,Almog-14} precisely yield~\eqref{0-10} when combined with elementary homogenization theory.
The proof is based on (scalar) potential theory and crucially relies on the facts that the space dimension is $d=3$, that~$A$ is everywhere isotropic, and that the inclusions are spherical and disjoint.
Another contribution is due to Mourrat~\cite{Mourrat-13}, who studied for all~$d\ge 2$ a discrete elliptic equation (instead of a continuum one) with sparse i.i.d.\@ perturbations of the conductivity, proving~\eqref{0-10} in that setting by strongly relying on quantitative stochastic homogenization results of~\cite{MR2789576,MR2932541}. We also highlight the inspiring work~\cite{Anantharaman-LeBris-11,Anantharaman-LeBris-12} by Anantharaman and Le Bris, who obtained related results on sparse i.i.d.\@ perturbations of a periodic array of inclusions; see also~\cite{Anantharaman-these}.
Those different previous results were largely surpassed in~\cite{MR3458165}, where we studied the case of a general stationary and ergodic inclusion process, focusing on a dilute regime obtained by a Bernoulli deletion procedure where each inclusion is preserved independently with low probability, and where we established real analyticity with respect to the Bernoulli parameter --- only assuming that the number of inclusions per unit volume be uniformly bounded. This was recently extended in~\cite{doi:10.1142/S0219199722500274,CRMATH_2022__360_G8_909_0} to prove Gevrey regularity when starting from a Poisson point process (for which the uniform boundedness assumption fails).
In a different vein, partly inspired by Berdichevski{\u\i}'s approach in the periodic setting, Pertinand addressed in~\cite{Pertinand} the case when the dilution of the point process is obtained by dilating a given hardcore process, and he proved the real analyticity with respect to the inverse of the dilation parameter.
Very recently, G\'erard-Varet~\cite{MR4400909} proposed an alternative approach where
he bypasses homogenization theory and directly quantifies in terms of the volume fraction $\varphi$ the distance between the solution of a problem with sparse inclusions and that of an effective problem with conductivity given by the Clausius--Mossotti formula.
Various related contributions concern the validity of Einstein's formula for the effective viscosity of a Stokes fluid with a dilute suspension of rigid particles; see~\cite{MR2982744,MR4098775,MR4102716,MR4260456,2008.03837}.
Our approach in the present note is an adaptation of the recent short proof that we obtained in~\cite{2008.03837} for Einstein's formula: it allows to justify the Clausius--Mossotti formula for the first time in the setting of general point processes and to determine the optimal error estimate. Note that several points of the proof simplify in the present setting. In particular we manage to fully bypass the variational formulation of~\cite{2008.03837}.

While the Clausius--Mossotti formula is universal in the sense that it only depends on the set of inclusions via its volume fraction (and on the shape of the inclusions, here assumed to be spherical),
the next-order correction further depends on the two-point correlation function. The identification of this correction was first discussed in~\cite{MR3458165}, and it has been the object of many recent contributions in the context of Einstein's formula~\cite{BG-72,GVH,GV-20,GVM-20,2008.03837}; we refer in particular to our recent work~\cite{2008.03837} where all higher-order corrections are systematically described in form of a cluster expansion.

\section{Proof of Theorem~\ref{th:CM-DG+}}
We denote by $\Pc=\{x_n\}_n$ the point set sampled according to the underlying random point process.
In order to justify the approximation~\eqref{0-13} of the corrector in terms of single-inclusion problems, we start by singling out clusters of close inclusions:
let $\{K_{q}\}_q$ be the family of connected components of the fattened set $\Bc(\Pc)+ B$, and consider the corresponding clusters
\[\qquad J_{q}\,:=\,\bigcup_{n:B(x_n)\subset K_{q}}B(x_n).\]
As we focus here on spherical inclusions, we note that a point $x_n\in\Pc$ does not belong to a cluster $J_q$ made of at least two inclusions if and only if
\[\rho_n\,:=\,\tfrac12\inf_{m:m\ne n}|x_n-x_m| \,\ge\,2.\]
Let then $\Sc:=\{n:\rho_n \ge2\}$ be the set of indices corresponding to well-separated inclusions
and let~$\Sc'$ be the set of indices $q$ such that the cluster $J_{q}$ is made of at least two inclusions. (In what follows,~$\Sc$ or $\Sc'$ can be empty.)
Next, in order to define suitable neighborhoods of the inclusions, we consider the Voronoi tessellation~$\{V_n\}_{n}$ associated with the point set~$\Pc=\{x_n\}_n$, that is,
\begin{equation}\label{eq:Voronoi}
V_n\,:=\,\Big\{z\in\R^d:|z-x_n|<\inf_{m:m\ne n}|z-x_m|\Big\}.
\end{equation}
We merge Voronoi cells that correspond to inclusions in the same cluster, setting
\[W_q\,:=\,\bigcup_{n:B(x_n)\subset J_{q}}V_n,\]
and we then partition the whole space as
\[\R^d\,=\,\Big(\bigcup_{n\in\Sc}V_n\Big)\cup\Big(\bigcup_{q\in\Sc'}W_q\Big).\]
We shall repeatedly use the following elementary property of Voronoi tessellations, see~\cite[Lemma~2.5]{2008.03837}: for all stationary random fields $\zeta$ with $\expec{|\zeta|}<\infty$, we have
\begin{equation}\label{eq:Voronoi}
\expec{\zeta}\,=\,\expecM{\sum_{n\in\Sc}\frac{\mathds1_{0\in B(x_n)}}{|B|}\int_{V_n}\zeta}+\expecM{\sum_{q\in\Sc'}\frac{\mathds1_{0\in J_q}}{|J_q|}\int_{W_q}\zeta}.
\end{equation}
We shall also use the short-hand notation $\tilde\Bc(\Pc):=\cup_{n\in \Sc} B(x_n)$ for the union of well-separated inclusions.
We split the proof into six steps. Let $e\in\R^d$ with $|e|=1$.

\medskip
\step1 Representation formula for the error: proof that
\begin{equation}\label{eq:error-estim0}
\big|\bar Ae-(A_1e+\varphi\widehat A_2e)\big|
\,\lesssim\,
\expecM{\sum_{n\in\Sc}\mathds1_{B(x_n)}|\nabla(\phi_e-\psi_{e,n})|}
+\expecm{\mathds1_{\Bc(\Pc)\setminus\tilde\Bc(\Pc)}\big(1+|\nabla \phi_e|\big)},
\end{equation}
where $\widehat A_2$ is given by~\eqref{eq:def-corrA2} in the statement, and where for all $n$ we define $\psi_{e,n}$ as the unique weak solution in $\dot H^1(\R^d)$ of the whole-space single-inclusion problem centered at~$x_n$, that is, 
\begin{equation*} 
-\nabla \cdot \big(A_1+\mathds1_{B(x_n)}(A_2-A_1)\big)(\nabla \psi_{e,n}+e)=0,\qquad\text{in $\R^d$}.
\end{equation*}
Recall that as in the statement we also denote by $\psi_e$ the solution of the corresponding single-inclusion problem with center~$x_n$ replaced by the origin~$0$, cf.~\eqref{eq:single-incl-1}.

\medskip\noindent
By definition of the two-phase coefficient field $A$, cf.~\eqref{eq:model-A}, and of the associated effective conductivity~$\bar A$, cf.~\eqref{0-1},  we find
\begin{eqnarray*}
\bar A e& =& \expec{A(\nabla \phi_e+e)}\\
&=&\expec{A_1 (\nabla \phi_e+e)}+\expec{\mathds1_{\Bc(\Pc)}(A_2-A_1) (\nabla \phi_e+e)}\\
&=& A_1e+\expec{\mathds1_{\Bc(\Pc)}(A_2-A_1) (\nabla \phi_e+e)},
\end{eqnarray*}
where the last identity follows from the fact that $A_1$ is constant and $\expec{\nabla\phi_e}=0$.
Decomposing $\Bc(\Pc)$ into well-separated inclusions and into clusters of close inclusions, and comparing the corrector $\phi_e$ to single-inclusion solutions $\{\psi_{e,n}\}_n$ in well-separated inclusions, we can decompose
\begin{multline}\label{eq:pre-decomp-barA}
\bar A e
\,=\, A_1e+\expecM{\sum_{n\in\Sc}\mathds1_{B(x_n)}(A_2-A_1) (\nabla\psi_{e,n}+e)}\\
+\expecM{\sum_{n\in\Sc}\mathds1_{B(x_n)}(A_2-A_1) \nabla(\phi_e-\psi_{e,n})}
+\expec{\mathds1_{\Bc(\Pc)\setminus\tilde\Bc(\Pc)}(A_2-A_1) (\nabla \phi_e+e)}.
\end{multline}
In order to reformulate the first right-hand side term, we appeal to~\eqref{eq:Voronoi}, to the effect of
\begin{equation*}
\expecM{\sum_{n\in\Sc}\mathds1_{B(x_n)}(A_2-A_1) (\nabla\psi_{e,n}+e)}
\,=\,\expecM{\sum_{n\in\Sc}\frac{\mathds1_{0\in B(x_n)}}{|B|}\int_{B(x_n)}(A_2-A_1) (\nabla\psi_{e,n}+e)},
\end{equation*}
which can be further rewritten as follows, recalling that $A_2$ is stationary and independent of the point process,
\begin{eqnarray*}
\expecM{\sum_{n\in\Sc}\mathds1_{B(x_n)}(A_2-A_1) (\nabla\psi_{e,n}+e)}
&=&\expecM{\sum_{n\in\Sc}\frac{\mathds1_{0\in B(x_n)}}{|B|}}\expecM{\int_{B}(A_2-A_1) (\nabla\psi_{e}+e)}\\
&=&\expecm{\mathds1_{\tilde\Bc(\Pc)}}\expecM{\fint_{B}(A_2-A_1) (\nabla\psi_{e}+e)}.
\end{eqnarray*}
Recognizing the definition~\eqref{eq:def-corrA2} of~$\widehat A_2$, noting that the definition~\eqref{eq:def-varphi} of the volume fraction yields
\[\expecm{\mathds1_{\tilde\Bc(\Pc)}}\,=\,\varphi-\expecm{\mathds1_{\Bc(\Pc)\setminus\tilde\Bc(\Pc)}},\]
and using the energy estimate $\int_{\R^d}|\nabla\psi_e|^2\lesssim1$ for the solution of~\eqref{eq:single-incl-1},
we deduce
\begin{eqnarray*}
\bigg|\expecM{\sum_{n\in\Sc}\mathds1_{B(x_n)}(A_2-A_1) (\nabla\psi_{e,n}+e)}-\varphi\widehat A_2e\bigg|
\,\lesssim\,\expecm{\mathds1_{\Bc(\Pc)\setminus\tilde\Bc(\Pc)}}.
\end{eqnarray*}
Combined with~\eqref{eq:pre-decomp-barA}, this yields the claim~\eqref{eq:error-estim0}.

\medskip
\step2 Approximation of the corrector by local Dirichlet problems.\\
Instead of directly comparing the corrector $\phi_e$ to whole-space single-inclusion solutions $\{\psi_{e,n}\}_n$ as in~\eqref{eq:error-estim0}, we start by comparing to solutions of single-inclusion Dirichlet problems in Voronoi cells.
More precisely, for all $n\in\Sc$, we define $\psi_{e,n}^\circ$ as the unique weak solution in $H^1_0(V_n)$ of the single-inclusion problem
\begin{equation}\label{eq:psi0n}
-\nabla \cdot A (\nabla \psi_{e,n}^\circ+e)=0,\qquad\text{in $V_n$},
\end{equation}
and, for all $q\in\Sc'$, we define $\psi_{e,q}^{\prime\circ}$ as the unique weak solution in $H^1_0(W_q)$ of
\[-\nabla \cdot A (\nabla \psi_{e,q}^{\prime\circ}+e)=0,\qquad\text{in $W_q$}.\]
Implicitly extending $\psi_{e,n}^\circ$ by zero outside $V_n$, and $\psi^{\prime\circ}_{e,q}$ by zero outside $W_q$, we then set
\[\psi_e^\circ\,:=\,\sum_{n\in \Sc} \psi_{e,n}^\circ+\sum_{q \in W_q} \psi_{e,q}^{\prime\circ}.\]
By definition, $\nabla \psi_e^\circ$ is stationary and we claim that it has vanishing expectation and finite second moments,
\begin{equation}\label{eq:prop-nab-psi0}
\expec{\nabla\psi_e^\circ}=0,\qquad\expec{|\nabla \psi_e^\circ|^2} \lesssim 1.
\end{equation}
Indeed, for all $R>0$, applying~\eqref{eq:Voronoi} to $|\nabla\psi_e^\circ|^2\wedge R$, we find
\begin{equation*}
\expec{|\nabla \psi_e^\circ|^2\wedge R}
\,\le\,\expecM{\sum_{n\in \Sc} \frac{\mathds{1}_{0\in B(x_n)}}{|B|} \int_{V_n} |\nabla \psi_e^\circ|^2}+ \expecM{\sum_{q\in \Sc'} \frac{\mathds{1}_{0\in J_{q}}}{|J_{q}|} \int_{W_q} |\nabla \psi_e^\circ|^2},
\end{equation*}
and thus, using energy estimates for $\psi_e^\circ$ in Voronoi cells, and further applying~\eqref{eq:Voronoi} to the constant function $1$,
\begin{equation*}
\expec{|\nabla \psi_e^\circ|^2\wedge R}
\,\lesssim\,\expecM{\sum_{n\in \Sc} \frac{\mathds{1}_{0\in B(x_n)}}{|B|} |V_n|}+ \expecM{\sum_{q\in \Sc'} \frac{\mathds{1}_{0\in J_{q}}}{|J_{q}|} |W_q|}
\,=\,1.
\end{equation*}
By the monotone convergence theorem, this proves the claim $\expec{|\nabla\psi_e^\circ|^2}\lesssim1$.
Next, we can apply~\eqref{eq:Voronoi} to $\nabla\psi_e^\circ$, to the effect of
\begin{eqnarray*}
\expec{\nabla \psi_e^\circ}
&=& \expecM{\sum_{n\in \Sc} \frac{\mathds{1}_{0\in B(x_n)}}{|B|} \int_{V_n} \nabla \psi_e^\circ}+ \expecM{\sum_{q\in \Sc'} \frac{\mathds{1}_{0\in J_{q}}}{|J_{q}|} \int_{W_q} \nabla \psi_e^\circ},
\end{eqnarray*}
where the right-hand side vanishes due to homogeneous Dirichlet boundary conditions, and the claim~\eqref{eq:prop-nab-psi0} follows.

\medskip
\step3 Approximation error estimate: proof that
\begin{equation}\label{0-20}
\expec{ |\nabla (\phi_e-\psi_e^\circ)|^2}+\expecm{ \mathds1_{\Bc(\Pc)} |\nabla (\phi_e-\psi_e^\circ)|}
\,\lesssim\,\expecM{\sum_{n\in\Sc}\frac{\mathds1_{0\in B(x_n)}}{|B|}\rho_n^{-d}}+\expecm{\mathds1_{\Bc(\Pc)\setminus\tilde\Bc(\Pc)}}.
\end{equation}
We start by proving the estimate on $\expec{ |\nabla (\phi_e-\psi_e^\circ)|^2}$.
Using the corrector equation~\eqref{0-2} for~$\phi_e$ in form of
\[\expec{ \nabla (\phi_e-\psi_e^\circ)\cdot A(\nabla\phi_e+e)}=0,\]
we find
\begin{eqnarray*}
\expec{ |\nabla (\phi_e-\psi_e^\circ)|^2}
&\lesssim&\expec{ \nabla (\phi_e-\psi_e^\circ)\cdot A\nabla (\phi_e-\psi_e^\circ)}\\
&=&-\expec{ \nabla (\phi_e-\psi_e^\circ)\cdot A(\nabla\psi_e^\circ+e)},
\end{eqnarray*}
which we can further decompose into
\begin{multline*}
\expec{ |\nabla (\phi_e-\psi_e^\circ)|^2}
\,\lesssim\,
-\expecm{\nabla(\phi_e-\psi_e^\circ)\cdot Ae}\\
+\expecm{(\nabla\psi_e^\circ+e)\cdot A\nabla\psi_e^\circ}
-\expecm{(\nabla\phi_e+e)\cdot A\nabla\psi_e^\circ}.
\end{multline*}
As $A_1$ is constant and as $\expec{\nabla(\phi_e-\psi_e^\circ)}=0$, the first right-hand side term is equal to
\begin{equation*}
\expecm{\nabla(\phi_e-\psi_e^\circ)\cdot Ae}\,=\,\expecm{\nabla(\phi_e-\psi_e^\circ)\cdot \mathds1_{\Bc(\Pc)}(A_2-A_1)e},
\end{equation*}
hence
\begin{multline}\label{eq:pre-decomp-phi-psi0}
\expec{ |\nabla (\phi_e-\psi_e^\circ)|^2}
\,\lesssim\,
-\expecm{\nabla(\phi_e-\psi_e^\circ)\cdot \mathds1_{\Bc(\Pc)}(A_2-A_1)e}\\
+\expecm{(\nabla\psi_e^\circ+e)\cdot A\nabla\psi_e^\circ}
-\expecm{(\nabla\phi_e+e)\cdot A\nabla\psi_e^\circ}.
\end{multline}
If $A$ was symmetric, then the corrector equation would ensure that the last right-hand side term vanishes. Using that $A_1$ is constant, we shall show that, even though this term does not vanish in the general non-symmetric case, it can be localized inside inclusions.
For that purpose, we rewrite the corrector equation as
\[-\nabla\cdot A_1\nabla\phi_e\,=\,\nabla\cdot Ae+\nabla\cdot \mathds1_{\Bc(\Pc)}(A_2-A_1)\nabla\phi_e,\qquad\text{in $\R^d$},\]
and we note that, as the coefficient $A_1$ is constant, it can be replaced by its pointwise transpose $A_1^T$ in the left-hand side,
\[-\nabla\cdot A_1^T\nabla\phi_e\,=\,\nabla\cdot Ae+\nabla\cdot \mathds1_{\Bc(\Pc)}(A_2-A_1)\nabla\phi_e\qquad\text{in $\R^d$}.\]
Testing this equation with $\psi_e^\circ$ then yields
\[\expecm{\nabla\psi_e^\circ\cdot A_1^T\nabla\phi_e}\,=\,-\expecm{\nabla\psi_e^\circ\cdot Ae}-\expecm{\nabla\psi_e^\circ\cdot\mathds1_{\Bc(\Pc)}(A_2-A_1)\nabla\phi_e},\]
or equivalently, adding and subtracting several terms,
\begin{multline*}
\expecm{(\nabla\phi_e+e)\cdot A\nabla\psi_e^\circ}
\,=\,\expecm{e\cdot A\nabla\psi_e^\circ}
-\expecm{\nabla\psi_e^\circ\cdot Ae}\\
+\expecm{\nabla(\phi_e-\psi_e^\circ)\cdot \mathds1_{\Bc(\Pc)}(A_2-A_1)\nabla\psi_e^\circ}
-\expecm{\nabla\psi_e^\circ\cdot\mathds1_{\Bc(\Pc)}(A_2-A_1)\nabla(\phi_e-\psi_e^\circ)}.
\end{multline*}
Inserting this into~\eqref{eq:pre-decomp-phi-psi0}, we get
\begin{multline}\label{eq:pre-decomp-phi-psi1}
\expec{ |\nabla (\phi_e-\psi_e^\circ)|^2}
\,\lesssim\,
\expecm{\nabla\psi_e^\circ\cdot A(\nabla\psi_e^\circ+e)}
-\expecm{\nabla(\phi_e-\psi_e^\circ)\cdot \mathds1_{\Bc(\Pc)}(A_2-A_1)e}\\
+\expecm{\nabla\psi_e^\circ\cdot\mathds1_{\Bc(\Pc)}(A_2-A_1)\nabla(\phi_e-\psi_e^\circ)}\\
-\expecm{\nabla(\phi_e-\psi_e^\circ)\cdot \mathds1_{\Bc(\Pc)}(A_2-A_1)\nabla\psi_e^\circ}.
\end{multline}
We note that the first right-hand side term vanishes: indeed, appealing to~\eqref{eq:Voronoi}, we find
\begin{eqnarray*}
\lefteqn{\expecm{ \nabla\psi_e^\circ\cdot A(\nabla\psi_e^\circ+e)}}\\
&=&\expecM{\sum_{n\in\Sc}\frac{\mathds1_{0\in B(x_n)}}{|B|}\int_{V_n} \nabla\psi_e^\circ\cdot A(\nabla\psi_e^\circ+e)}+\expecM{\sum_{q\in\Sc'}\frac{\mathds1_{0\in J_q}}{|J_q|}\int_{W_q} \nabla\psi_e^\circ\cdot A(\nabla\psi_e^\circ+e)},
\end{eqnarray*}
where for all $n\in\Sc$ and $q\in\Sc'$ the equations for $\psi_e^\circ|_{V_n}=\psi_{e,n}^\circ$ and $\psi^\circ_e|_{W_q}=\psi_{e,q}^{\prime\circ}$ precisely give
\begin{eqnarray*}
\int_{V_n} \nabla\psi_e^\circ\cdot A(\nabla\psi_e^\circ+e)\,=\,0,\qquad\int_{W_q} \nabla\psi_e^\circ\cdot A(\nabla\psi_e^\circ+e)\,=\,0.
\end{eqnarray*}
The estimate~\eqref{eq:pre-decomp-phi-psi1} then leads us to
\[\expec{ |\nabla (\phi_e-\psi_e^\circ)|^2}
\,\lesssim\,\expecm{\mathds1_{\Bc(\Pc)}\big(1+|\nabla\psi_e^\circ|\big)|\nabla (\phi_e-\psi_e^\circ)|}.\]
Appealing to~\eqref{eq:Voronoi}, together with the Cauchy--Schwarz inequality, this yields
\begin{multline}\label{eq:decomp-phi-psi1}
\expec{ |\nabla (\phi_e-\psi_e^\circ)|^2}
\,\lesssim\,\expecM{\sum_{n\in\Sc}\frac{\mathds1_{0\in B(x_n)}}{|B|}\Big(1+\int_{B(x_n)}|\nabla\psi_e^\circ|^2\Big)^\frac12\Big(\int_{B(x_n)} |\nabla (\phi_e-\psi_e^\circ)|^2\Big)^\frac12}\\
+\expecM{\sum_{q\in\Sc'}\frac{\mathds1_{0\in J_q}}{|J_q|}\Big(|J_q|+\int_{J_q}|\nabla\psi_e^\circ|^2\Big)^\frac12\Big(\int_{J_q}|\nabla (\phi_e-\psi_e^\circ)|^2\Big)^\frac12}.
\end{multline}
For all $n\in\Sc$, as $A_1$ is constant, the defining equation~\eqref{eq:psi0n} for $\psi_e^\circ|_{V_n}=\psi_{e,n}^\circ\in H^1_0(V_n)$ yields
\begin{equation*}
\int_{V_n}|\nabla\psi_e^\circ|^2\,\lesssim\,\int_{V_n}\nabla\psi_{e,n}^\circ\cdot A\nabla\psi_{e,n}^\circ
\,=\,-\int_{V_n}\nabla\psi_{e,n}^\circ\cdot Ae
\,=\,-\int_{B(x_n)}\nabla\psi_{e,n}^\circ\cdot (A_2-A_1)e,
\end{equation*}
which implies the energy estimate
\[\int_{V_n}|\nabla\psi_e^\circ|^2\,\lesssim\,1.\]
Similarly, for all $q\in\Sc'$, we find
\begin{equation}\label{eq:energy-Wq-psi0}
\int_{W_q}|\nabla\psi_e^\circ|^2\,\lesssim\,|J_q|.
\end{equation}
Inserting this into~\eqref{eq:decomp-phi-psi1}, we deduce
\begin{multline*}
\expec{ |\nabla (\phi_e-\psi_e^\circ)|^2}
\,\lesssim\,\expecM{\sum_{n\in\Sc}\frac{\mathds1_{0\in B(x_n)}}{|B|}\Big(\int_{B(x_n)} |\nabla (\phi_e-\psi_e^\circ)|^2\Big)^\frac12}\\
+\expecM{\sum_{q\in\Sc'}\frac{\mathds1_{0\in J_q}}{|J_q|}|J_q|^\frac12\Big(\int_{W_q}|\nabla (\phi_e-\psi_e^\circ)|^2\Big)^\frac12}.
\end{multline*}
We now appeal to the following mean-value property, which we shall prove in Step~4 below,
\begin{equation}\label{0-16}
\int_{B(x_n)} |\nabla (\phi_e-\psi_e^\circ)|^2 \,\lesssim\, \rho_n^{-d} \int_{V_n}  |\nabla (\phi_e-\psi_e^\circ)|^2.
\end{equation}
By Young's inequality, we then get for any $K>0$,
\begin{multline*}
\expec{ |\nabla (\phi_e-\psi_e^\circ)|^2}
\,\lesssim\,K\expecM{\sum_{n\in\Sc}\frac{\mathds1_{0\in B(x_n)}}{|B|}\rho_n^{-d}}+K\expecM{\sum_{q\in\Sc'}\frac{\mathds1_{0\in J_q}}{|J_q|}|J_q|}\\
+\tfrac1K\expecM{\sum_{n\in\Sc}\frac{\mathds1_{0\in B(x_n)}}{|B|}\int_{V_n}|\nabla (\phi_e-\psi_e^\circ)|^2}+\tfrac1K\expecM{\sum_{q\in\Sc'}\frac{\mathds1_{0\in J_q}}{|J_q|}\int_{W_q}|\nabla (\phi_e-\psi_e^\circ)|^2}.
\end{multline*}
As~\eqref{eq:Voronoi} implies that the last two right-hand side terms are equal to $\frac1K\expec{ |\nabla (\phi_e-\psi_e^\circ)|^2}$, choosing $K\simeq1$ large enough, and recalling the notation $\Bc(\Pc)\setminus\tilde\Bc(\Pc)=\cup_{q\in\Sc'}J_q$, we finally obtain
\begin{equation}\label{0-20/a}
\expec{ |\nabla (\phi_e-\psi_e^\circ)|^2}\,\lesssim\,\expecM{\sum_{n\in\Sc}\frac{\mathds1_{0\in B(x_n)}}{|B|}\rho_n^{-d}}+\expecm{\mathds1_{\Bc(\Pc)\setminus\tilde\Bc(\Pc)}}.
\end{equation}
To conclude the proof of the claim~\eqref{0-20}, it remains to establish the corresponding estimate on $\expec{\mathds1_{\Bc(\Pc)}|\nabla (\phi_e-\psi_e^\circ)|}$.
For that purpose, we start by appealing to~\eqref{eq:Voronoi} in form of
\begin{multline*}
\expec{\mathds1_{\Bc} |\nabla (\phi_e-\psi^\circ)|}
\\
\,=\,\expecM{\sum_{n\in \Sc} \frac{\mathds{1}_{0 \in B(x_n)}}{|B|} \int_{B(x_n)}  |\nabla (\phi_e-\psi^\circ)| }+ \expecM{\sum_{q\in \Sc'} \frac{\mathds{1}_{0 \in J_{q}}}{|J_{q}|} \int_{J_{q}} |\nabla(\phi_e-\psi^\circ)|}.
\end{multline*}
By the Cauchy--Schwarz inequality and the mean-value property~\eqref{0-16}, we find
\begin{multline*}
\expec{\mathds1_{\Bc} |\nabla (\phi_e-\psi_e^\circ)|}
\,\lesssim\,\expecM{\sum_{n\in \Sc} \frac{\mathds{1}_{0 \in B(x_n)}}{|B|}\rho_n^{-\frac d2}\Big(\int_{V_n}  |\nabla (\phi_e-\psi_e^\circ)|^2\Big)^\frac12}\\
+\expecM{\sum_{q\in\Sc'}\frac{\mathds1_{0\in J_q}}{|J_q|}|J_q|^\frac12\Big(\int_{W_q}|\nabla (\phi_e-\psi_e^\circ)|^2\Big)^\frac12}.
\end{multline*}
By Young's inequality, arguing as above, we obtain
\begin{equation*}
\expec{\mathds1_{\Bc} |\nabla (\phi_e-\psi_e^\circ)|}
\,\lesssim\,\expecM{\sum_{n\in \Sc} \frac{\mathds{1}_{0 \in B(x_n)}}{|B|}\rho_n^{-d}}+\expecm{\mathds1_{\Bc(\Pc)\setminus\tilde\Bc(\Pc)}}+\expec{|\nabla (\phi_e-\psi_e^\circ)|^2},
\end{equation*}
and the claim~\eqref{0-20} now follows from~\eqref{0-20/a}.

\medskip
\step4 Proof of the mean-value property~\eqref{0-16}. \\
We reformulate \eqref{0-16} in the following form: for all $r\ge1$, for all $v\in H^1(rB)$ that satisfies in the weak sense,
\[-\nabla \cdot \big(A_1+\mathds 1_{B}(A_2-A_1)\big) \nabla v=0,\qquad\text{in $rB$},\]
we have
\begin{equation}\label{0-21} 
\int_{B} |\nabla v|^2\, \lesssim\, r^{-d} \int_{rB} |\nabla v|^2.
\end{equation}
To prove this result, we decompose $v=v_1+v_2$, where $v_1$ is the unique weak solution in $v+H^1_0(rB)$ of
\[-\nabla \cdot A_1 \nabla v_1=0,\qquad\text{in $rB$},\]
and where $v_2$ is the unique weak solution in $H^1_0(rB)$ of
$$
-\nabla \cdot\big(A_1+\mathds 1_{B}(A_2-A_1)\big)  \nabla v_2 =\nabla \cdot\mathds1_B(A_2-A_1)\nabla v_1,\qquad\text{in $rB$}.
$$
As $A_1$ is constant, we can apply to $v_1$ the mean-value property for harmonic functions, to the effect of
\[\int_{B} |\nabla v_1|^2\, \lesssim\, r^{-d} \int_{rB} |\nabla v_1|^2.\]
As the defining Dirichlet problem for $v_1$ yields the energy estimate
\[\int_{rB}|\nabla v_1|^2\lesssim\int_{rB}|\nabla v|^2,\]
we deduce
\begin{equation}\label{0-22a}
\int_{B} |\nabla v_1|^2\, \lesssim\, r^{-d} \int_{rB} |\nabla v|^2.
\end{equation}
In order to estimate $v_2$, we start from the corresponding energy estimate
\[\int_{rB} |\nabla v_2|^2 \, \lesssim\, \int_{B} |\nabla v_1|^2.\]
Combined with~\eqref{0-22a}, this yields the claim~\eqref{0-21} by the triangle inequality.

\medskip
\step5 Comparison of whole-space and of Dirichlet single-inclusion problems: proof that for all $n\in \Sc$ we have
\begin{equation}\label{0-22}
\int_{B(x_n)} |\nabla (\psi_{e,n}-\psi^\circ_{e})|^2 \,\lesssim\, \rho_n^{-2d}.
\end{equation}
By definition, the difference $\delta \psi_{e,n}:=\psi_{e,n}-\psi^\circ_{e,n}$ satisfies $\delta\psi_{e,n}|_{V_n}\in\psi_{e,n}+H^1_0(V_n)$ and, in the weak sense,
\begin{equation}\label{eq:delta-psin}
-\nabla  \cdot A \nabla \delta \psi_{e,n}=0,\qquad\text{in $V_n$}.
\end{equation}
The mean-value property~\eqref{0-21} on $B_{\rho_n}(x_n)\subset V_n$ then yields
\begin{equation}\label{eq:estim-0st4}
\int_{B(x_n)} |\nabla \delta \psi_{e,n}|^2 \,\lesssim \,  \rho_n^{-d} \int_{B_{\rho_n}(x_n)} |\nabla \delta \psi_{e,n}|^2\,\le\,\rho_n^{-d} \int_{V_n} |\nabla \delta \psi_{e,n}|^2.
\end{equation}
Given a smooth cut-off $\chi_n$ such that
\[\chi_n|_{B_{\frac12\rho_n}(x_n)}\,=\,0,\qquad\chi_n|_{\R^d\setminus B_{\rho_n}(x_n)}\,=\,1,\qquad|\nabla \chi_n|\,\lesssim\, \rho_n^{-1},\]
the Dirichlet problem~\eqref{eq:delta-psin} for $\delta\psi_{e,n}$ implies the energy estimate
\[\int_{V_n}|\nabla\delta \psi_{e,n}|^2\,\lesssim\,\int_{V_n}|\nabla(\chi_n \psi_{e,n})|^2,\]
and thus, as the single-inclusion solution $\psi_{e,n}$ enjoys the same decay properties as in~\eqref{0-11}, we deduce
\[\int_{V_n}  |\nabla \delta \psi_{e,n}|^2\,\lesssim\,  \rho_n^{-d}.\]
Combined with~\eqref{eq:estim-0st4}, this proves the claim~\eqref{0-22}.

\medskip
\step6 Conclusion.\\
Starting from~\eqref{eq:error-estim0} and comparing single-inclusion solutions to their Dirichlet version $\psi_e^\circ$, we have
\begin{multline*}
\big|\bar Ae-(A_1e+\varphi\widehat A_2e)\big|
\,\lesssim\,
\expecm{\mathds1_{\Bc(\Pc)}|\nabla(\phi_e-\psi_e^\circ)|}
+\expecm{\mathds1_{\Bc(\Pc)\setminus\tilde\Bc(\Pc)}\big(1+|\nabla \psi^\circ_e|\big)}\\
+\expecM{\sum_{n\in\Sc}\mathds1_{B(x_n)}|\nabla(\psi_{e,n}-\psi_e^\circ)|}.
\end{multline*}
Using~\eqref{0-20} to estimate the first right-hand side term, appealing to~\eqref{eq:Voronoi} together with the energy estimate~\eqref{eq:energy-Wq-psi0} in form of
\begin{equation*}
\expecm{\mathds1_{\Bc(\Pc)\setminus\tilde\Bc(\Pc)}|\nabla \psi_e^\circ|}\,=\,\expecM{\sum_{q\in\Sc'}\frac{\mathds1_{0\in J_q}}{|J_q|}\int_{J_q}|\nabla \psi_e^\circ|}
\,\lesssim\,\expecM{\sum_{q\in\Sc'}\frac{\mathds1_{0\in J_q}}{|J_q|}|J_q|}\,\le\,\expecm{\mathds1_{\Bc(\Pc)\setminus\tilde\Bc(\Pc)}},
\end{equation*}
and appealing to~\eqref{eq:Voronoi} together with~\eqref{0-22} in form of
\begin{eqnarray*}
\expecM{\sum_{n\in\Sc}\mathds1_{B(x_n)}|\nabla (\psi_{e,n}-\psi_{e}^\circ)|}&=&\expecM{\sum_{n\in\Sc}\frac{\mathds1_{0\in B(x_n)}}{|B|}\int_{B(x_n)}|\nabla (\psi_{e,n}-\psi_{e}^\circ)|}\\
&\lesssim&\expecM{\sum_{n\in\Sc}\frac{\mathds1_{0\in B(x_n)}}{|B|}\rho_n^{-d}},
\end{eqnarray*}
we deduce
\begin{equation}\label{bingo}
\big|\bar A e- (A_1e+\varphi \widehat A_2e)\big|
\,\lesssim\,
\expecM{\sum_{n}\frac{\mathds1_{0\in B(x_n)}}{|B|}\langle\rho_n\rangle^{-d}}
+\expecm{\mathds1_{\Bc(\Pc)\setminus\tilde\Bc(\Pc)}}.
\end{equation}
It remains to evaluate the two right-hand side terms.
For that purpose, we shall appeal to the following observation that we first made in~\cite{2008.03837}: for any non-increasing function \mbox{$g\in\Ld^\infty(\R^+)$} with $g(r)\downarrow0$ as $r\uparrow\infty$, there holds
\begin{equation}\label{0-23}
\expecM{\sum_{n} \mathds1_{0 \in B(x_n)} g(\rho_n)}\,\lesssim\, \lambda_2 \|g\|_{\Ld^\infty(\R^+)}+\int_{0}^\infty |g'(r)|\,\Big((\lambda_2 \,\langle r\rangle^d)\wedge \lambda \Big)\,dr,
\end{equation}
with the short-hand notation $\langle r\rangle:=(1+r^2)^\frac12$.
Applying this to $g(r)=\langle r\rangle^{-d}$, we find
\begin{eqnarray*}
\expecM{\sum_{n} \mathds{1}_{0 \in B(x_n)} \langle\rho_n\rangle^{-d}}&\lesssim&\lambda_2+\int_{0}^\infty\langle r\rangle^{-d-1}\Big((\lambda_2\langle r\rangle^d)\wedge \lambda \Big)\,dr\\
&\lesssim &
\lambda_2 \log(2+\tfrac{\lambda}{\lambda_2}) ~\lesssim~ \lambda_2 |\!\log \lambda_2|.
\end{eqnarray*}
Applying it to $g(r)=\mathds1_{r\le2}$, we further get
\[\expecm{\mathds 1_{\Bc(\Pc)\setminus\tilde \Bc(\Pc)}} \,=\,\expecM{\sum_n\mathds1_{0\in B(x_n)}\mathds1_{\rho_n\le2}}\,\lesssim\,
 \lambda_2.\]
Combined with~\eqref{bingo}, this concludes the proof of~\eqref{0-26+}.

\medskip\noindent
Finally, for completeness, we include a short proof of~\eqref{0-23}. For that purpose, we start by rewriting the left-hand side as 
\begin{equation}\label{0-24}
\expecM{\sum_{n}\mathds1_{0\in B(x_n)} g(\rho_n)}\,=\,\int_0^\infty g(r)\,d \Lambda(r),
\end{equation}
where the positive measure $\Lambda$ on $\R^+$ is defined by its distribution function
\[\Lambda([0,r])\,:=\,\expecM{\sum_{n} \mathds1_{0\in B(x_n)} \mathds1_{\rho_n\le r}}
\,=\,\expecM{\sum_{n}\mathds1_{|x_n|<1}\,\mathds1_{\exists m\ne n:\,|x_m-x_n|\le 2r}}.\]
By definition of the minimal length $\ell$ of the point process, cf.~\eqref{0-14}, note that $\Lambda([0,r])=0$ for $r<\frac12\ell$.
Moreover, we can bound
\[\Lambda([0,r])\,\le\, \expecM{\sum_{n} \mathds1_{|x_n|<1}}\,=\, \lambda|B|,\]
and alternatively, for $r\ge\frac12\ell$, by definition of the second-order intensity, cf.~\eqref{eq:high-intens/Re},
\begin{multline*}
\Lambda([0,r])\,\le\, \expecM{\sum_{n\ne m}\mathds1_{|x_n|<1}\,\mathds1_{|x_m-x_n|\le 2r}}
\,=\,\iint_{B\times B_{2r}}f_2(x,x+y)\,dxdy\\
\,=\,(2r)^{-d}\iint_{B_{2r}\times B_{2r}}f_2(x,x+y)\,dxdy
\,\lesssim\,\lambda_2r^{d}.
\end{multline*}
Combining these estimates yields for all $r\ge0$,
\begin{equation}\label{0-25}
\Lambda([0,r])\,\lesssim\,
(\lambda_2 \langle r\rangle^d)\wedge \lambda .
\end{equation}
Under our assumptions on $g$, an integration by parts yields
\[\int_0^\infty g(r)\,d\Lambda(r)\,=\,-g(0)\Lambda(\{0\})+\int_{0}^\infty |g'(r)|\,\Lambda([0,r])\,dr,\]
and the claim~\eqref{0-23} follows in combination with~\eqref{0-24} and~\eqref{0-25}.
\qed

\section*{Acknowledgements}
MD acknowledges financial support from F.R.S.-FNRS, and AG from the European Research Council (ERC) under the European Union's Horizon 2020 research and innovation programme (Grant Agreement n$^\circ$~864066).

\bibliographystyle{plain}
\bibliography{biblio}
\end{document}